\documentclass[12pt,british]{amsart}
\usepackage{amsmath}
\usepackage{libertine}
\usepackage[libertine]{newtxmath}
\usepackage{mathrsfs}
\usepackage{hyperref}
\usepackage{babel}
\usepackage{amssymb,latexsym}
\usepackage{cite,amsmath,amstext,amsbsy,bm,url} 
\usepackage{paralist}
\usepackage{datetime} 
\usepackage{colortbl}
\usepackage[dvipsnames]{xcolor}
\usepackage{mathtools}
\usepackage{euscript}
\usepackage[all]{xy}
\usepackage{graphicx} 
\usepackage[a4paper,text={140mm,227mm},centering]{geometry} 
\usepackage{comment}  
\usepackage{cite}  
\usepackage{thmtools}
\usepackage{enumitem}
\usepackage{letltxmacro} 
\usepackage{nameref}
\usepackage{calc}
\usepackage{interval}
\usepackage{manfnt}
\usepackage{tikz-cd}  
\usepackage[utf8]{inputenc}
\usepackage{newunicodechar} 

\setdescription{labelindent=\parindent, leftmargin=2\parindent}
\setitemize[1]{labelindent=\parindent, leftmargin=2\parindent}
\setenumerate[1]{labelindent=0cm, leftmargin=*, widest=iiii}

\newunicodechar{α}{\ensuremath{\alpha}}
\newunicodechar{β}{\ensuremath{\beta}}
\newunicodechar{χ}{\ensuremath{\chi}}
\newunicodechar{δ}{\ensuremath{\delta}}
\newunicodechar{ε}{\ensuremath{\varepsilon}}
\newunicodechar{Δ}{\ensuremath{\Delta}}
\newunicodechar{η}{\ensuremath{\eta}}
\newunicodechar{γ}{\ensuremath{\gamma}}
\newunicodechar{Γ}{\ensuremath{\Gamma}}
\newunicodechar{ι}{\ensuremath{\iota}}
\newunicodechar{κ}{\ensuremath{\kappa}}
\newunicodechar{λ}{\ensuremath{\lambda}}
\newunicodechar{Λ}{\ensuremath{\Lambda}}
\newunicodechar{ν}{\ensuremath{\nu}}
\newunicodechar{μ}{\ensuremath{\mu}}
\newunicodechar{ω}{\ensuremath{\omega}}
\newunicodechar{Ω}{\ensuremath{\Omega}}
\newunicodechar{π}{\ensuremath{\pi}}
\newunicodechar{Π}{\ensuremath{\Pi}}
\newunicodechar{φ}{\ensuremath{\phi}}
\newunicodechar{Φ}{\ensuremath{\Phi}}
\newunicodechar{ψ}{\ensuremath{\psi}}
\newunicodechar{Ψ}{\ensuremath{\Psi}}
\newunicodechar{ρ}{\ensuremath{\rho}}
\newunicodechar{σ}{\ensuremath{\sigma}}
\newunicodechar{Σ}{\ensuremath{\Sigma}}
\newunicodechar{τ}{\ensuremath{\tau}}
\newunicodechar{θ}{\ensuremath{\theta}}
\newunicodechar{Θ}{\ensuremath{\Theta}}
\newunicodechar{ξ}{\ensuremath{\xi}}
\newunicodechar{Ξ}{\ensuremath{\Xi}}
\newunicodechar{ζ}{\ensuremath{\zeta}}

\newunicodechar{ℓ}{\ensuremath{\ell}}

\newunicodechar{𝔹}{\ensuremath{\bB}}
\newunicodechar{ℂ}{\ensuremath{\bC}}
\newunicodechar{𝔻}{\ensuremath{\bD}}
\newunicodechar{𝔼}{\ensuremath{\bE}}
\newunicodechar{𝔽}{\ensuremath{\bF}}
\newunicodechar{ℕ}{\ensuremath{\bN}}
\newunicodechar{ℙ}{\ensuremath{\bP}}
\newunicodechar{ℚ}{\ensuremath{\bQ}}
\newunicodechar{ℝ}{\ensuremath{\bR}}
\newunicodechar{𝕏}{\ensuremath{\bX}}
\newunicodechar{ℤ}{\ensuremath{\bZ}}
\newunicodechar{𝒜}{\ensuremath{\sA}}
\newunicodechar{ℬ}{\ensuremath{\sB}}
\newunicodechar{𝒞}{\ensuremath{\sC}}
\newunicodechar{𝒟}{\ensuremath{\sD}}
\newunicodechar{ℰ}{\ensuremath{\sE}}
\newunicodechar{ℱ}{\ensuremath{\sF}}
\newunicodechar{𝒢}{\ensuremath{\sG}}
\newunicodechar{ℋ}{\ensuremath{\sH}}
\newunicodechar{𝒥}{\ensuremath{\sJ}}
\newunicodechar{ℒ}{\ensuremath{\sL}}
\newunicodechar{𝒪}{\ensuremath{\sO}}
\newunicodechar{𝒬}{\ensuremath{\sQ}}
\newunicodechar{𝒯}{\ensuremath{\ts}}
\newunicodechar{𝒲}{\ensuremath{\sW}}

\newunicodechar{∂}{\ensuremath{\partial}}
\newunicodechar{∇}{\ensuremath{\nabla}}

\newunicodechar{↺}{\ensuremath{\circlearrowleft}}
\newunicodechar{∞}{\ensuremath{\infty}}
\newunicodechar{⊕}{\ensuremath{\oplus}}
\newunicodechar{⊗}{\ensuremath{\otimes}}
\newunicodechar{•}{\ensuremath{\bullet}}
\newunicodechar{Λ}{\ensuremath{\wedge}}
\newunicodechar{↪}{\ensuremath{\into}}
\newunicodechar{→}{\ensuremath{\to}}
\newunicodechar{↦}{\ensuremath{\mapsto}}
\newunicodechar{⨯}{\ensuremath{\times}}
\newunicodechar{∪}{\ensuremath{\cup}}
\newunicodechar{∩}{\ensuremath{\cap}}
\newunicodechar{⊋}{\ensuremath{\supsetneq}}
\newunicodechar{⊇}{\ensuremath{\supseteq}}
\newunicodechar{⊃}{\ensuremath{\supset}}
\newunicodechar{⊊}{\ensuremath{\subsetneq}}
\newunicodechar{⊆}{\ensuremath{\subseteq}}
\newunicodechar{⊂}{\ensuremath{\subset}}
\newunicodechar{≥}{\ensuremath{\geq}}
\newunicodechar{≠}{\ensuremath{\neq}}
\newunicodechar{≫}{\ensuremath{\gg}}
\newunicodechar{≪}{\ensuremath{\ll}}

\newunicodechar{≤}{\ensuremath{\leq}}
\newunicodechar{∈}{\ensuremath{\in}}
\newunicodechar{∉}{\ensuremath{\not \in}}
\newunicodechar{∖}{\ensuremath{\setminus}}
\newunicodechar{◦}{\ensuremath{\circ}}
\newunicodechar{°}{\ensuremath{^\circ}}
\newunicodechar{…}{\ifmmode\mathellipsis\else\textellipsis\fi}
\newunicodechar{·}{\ensuremath{\cdot}}
\newunicodechar{⋯}{\ensuremath{\cdots}}
\newunicodechar{∅}{\ensuremath{\emptyset}}
\newunicodechar{⇒}{\ensuremath{\Rightarrow}}

\newunicodechar{⁰}{\ensuremath{^0}}
\newunicodechar{¹}{\ensuremath{^1}}
\newunicodechar{²}{\ensuremath{^2}}
\newunicodechar{³}{\ensuremath{^3}}
\newunicodechar{⁴}{\ensuremath{^4}}
\newunicodechar{⁵}{\ensuremath{^5}}
\newunicodechar{⁶}{\ensuremath{^6}}
\newunicodechar{⁷}{\ensuremath{^7}}
\newunicodechar{⁸}{\ensuremath{^8}}
\newunicodechar{⁹}{\ensuremath{^9}}
\newunicodechar{ⁱ}{\ensuremath{^i}}

\newunicodechar{⌈}{\ensuremath{\lceil}}
\newunicodechar{⌉}{\ensuremath{\rceil}}
\newunicodechar{⌊}{\ensuremath{\lfloor}}
\newunicodechar{⌋}{\ensuremath{\rfloor}}

\newunicodechar{≅}{\ensuremath{\cong}}
\newunicodechar{⇔}{\ensuremath{\Leftrightarrow}}

\definecolor{linkred}{rgb}{0.7,0.2,0.2}
\definecolor{linkblue}{rgb}{0,0.2,0.6}

\usepackage[hyperpageref]{backref}

\DeclareFontFamily{OMS}{rsfs}{\skewchar\font'60}
\DeclareFontShape{OMS}{rsfs}{m}{n}{<-5>rsfs5 <5-7>rsfs7 <7->rsfs10 }{}
\DeclareSymbolFont{rsfs}{OMS}{rsfs}{m}{n}
\DeclareSymbolFontAlphabet{\scr}{rsfs}
\DeclareSymbolFontAlphabet{\scr}{rsfs}

\DeclareFontFamily{U}{mathx}{\hyphenchar\font45}
\DeclareFontShape{U}{mathx}{m}{n}{
	<5> <6> <7> <8> <9> <10>
	<10.95> <12> <14.4> <17.28> <20.74> <24.88>
	mathx10
}{}

\newtheorem{main}{Theorem}

 \newtheorem{theorem}{Theorem}[section]  
  \newtheorem{claim}{Claim}[theorem]
\newtheorem{conjecture}[theorem]{Conjecture}

\newtheorem{proposition}[theorem]{Proposition}
 \newtheorem{lemma}[theorem]{Lemma}
 
\theoremstyle{definition}
\newtheorem{definition}[theorem]{Definition}

\theoremstyle{remark}
\newtheorem{remark}[theorem]{Remark}

\newunicodechar{◦}{\ensuremath{\circ}}
\newunicodechar{°}{\ensuremath{^\circ}}

\newcommand{\sA}{\mathscr{A}}
\newcommand{\sB}{\mathscr{B}}
\newcommand{\sC}{\mathscr{C}}
\newcommand{\sD}{\mathscr{D}}
\newcommand{\sE}{\mathscr{E}}
\newcommand{\sF}{\mathscr{F}}
\newcommand{\sG}{\mathscr{G}}
\newcommand{\sH}{\mathscr{H}}

\newcommand{\sJ}{\mathscr{J}}

\newcommand{\sL}{\mathscr{L}}
\newcommand{\sM}{\mathscr{M}}

\newcommand{\sO}{\mathscr{O}}

\newcommand{\sQ}{\mathscr{Q}}

\newcommand{\sV}{\mathscr{V}}
\newcommand{\sW}{\mathscr{W}}
\newcommand{\sX}{\mathscr{X}}

\newcommand{\cV}{\mathcal V}

\newcommand{\kM}{\mathfrak M}

\newcommand{\bB}{\mathbb{B}}
\newcommand{\bC}{\mathbb{C}}
\newcommand{\bD}{\mathbb{D}}
\newcommand{\bE}{\mathbb{E}}
\newcommand{\bF}{\mathbb{F}}

\newcommand{\bN}{\mathbb{N}}

\newcommand{\bP}{\mathbb{P}}
\newcommand{\bQ}{\mathbb{Q}}
\newcommand{\bR}{\mathbb{R}}

\newcommand{\bV}{\mathbb{V}}

\newcommand{\bX}{\mathbb{X}}

\newcommand{\bZ}{\mathbb{Z}}

\makeatletter
\newcommand*{\rom}[1]{\expandafter\@slowromancap\romannumeral #1@}
\makeatother

\makeatletter
\hypersetup{
	pdfauthor={\authors},
	pdftitle={\@title},
	pdfsubject={\@subjclass},
	pdfkeywords={\@keywords},
	pdfstartview={Fit},
	pdfpagelayout={TwoColumnRight},
	pdfpagemode={UseOutlines},
	bookmarks,
	colorlinks,
	linkcolor=Sepia,
	citecolor=OliveGreen,
	urlcolor=WildStrawberry}
\makeatother

\begin{document} 
		\title[Hyperbolic   Campana's isotriviality conjecture]{Hyperbolicity of coarse moduli spaces\\ and isotriviality for certain families}

	\author{Ya Deng}
	\address{
		Department of Mathematical Sciences, Chalmers University of Technology and the University of Gothenburg, Sweden}
	\email{yade@chalmers.se\quad dengya.math@gmail.com}
	\urladdr{https://sites.google.com/site/dengyamath}
		\begin{abstract}
		In this paper, we prove the Kobayashi hyperbolicity of the coarse moduli spaces of canonically polarized or polarized Calabi-Yau manifolds in the sense of complex $V$-spaces (a generalization of   complex $V$-manifolds in the sense of Satake).  As an application, we prove the following hyperbolic version of Campana's isotriviality conjecture: for the smooth family of canonically polarized or polarized Calabi-Yau manifolds, when the Kobayashi pseudo-distance of the base vanishes identically, the family must be isotrivial, that is, any two fibers are isomorphic.   We also prove that for the  smooth projective family of polarized Calabi-Yau manifolds, its variation of the family is less than or equal to the essential dimension of the base. 
	\end{abstract}

	\maketitle
	
\tableofcontents
\section{Introduction}\label{introduction} 
\subsection{Main theorems}
Following Campana \cite[Definition 9.1]{Cam04}, a complex space $Y$ is called \emph{hyperbolically special} ({\it H}-special for short) if its Kobayashi pseudo distance $d_X:X\times X\to \mathbb{R}_{\geq 0}$  vanishes identically. The main motivation in this paper is to prove the following \emph{isotriviality} result for smooth families of canonically polarized manifolds or polarized Calabi-Yau manifolds.
 \begin{main}\label{main:isotrivial}
	Let $f:X\to Y $ a  proper smooth morphism of canonically polarized or polarized Calabi-Yau manifolds with connected fibers. Assume that  the base $Y$ is an {\it H}-special complex manifold. Then $f$ is \emph{isotrivial}, namely any two fibers of $f$ are isomorphic.
\end{main}

By the work of Fujiki-Schumacher \cite{FS90},  for moduli functor $\sM_H$ of canonically polarized manifolds or polarized Calabi-Yau manifolds with the   Hilbert polynomial $H(m)$, its coarse moduli space $\mathfrak{M}_H$ exists and is a \emph{Hausdorff reduced complex $V$-space}. The complex $V$-space generalizes the notation of $V$-manifold in the sense of Satake. As we will see in \S \ref{sec:V-hyper}, the definition of \emph{Kobayashi $V$-hyperbolicity} for complex $V$-manifolds introduced by Kobayashi \cite{Kob05} can also be generalized to complex $V$-spaces.

 With this notation in place, we state our more general result than Theorem \ref{main:isotrivial}.
\begin{main}\label{main}
The coarse moduli space $\mathfrak{M}_H$ of \emph{canonically polarized}   manifolds is \emph{Kobayashi $V$-hyperbolic}. In particular,	any $V$-morphism from an {\it H}-special manifold to $\mathfrak{M}_H$ is  necessarily constant.
\end{main}
The $V$-morphism from a complex manifold to the  complex $V$-space is, roughly speaking, a \emph{locally liftable} holomorphic map. As we shall see below, the \emph{Kobayashi pseudo distance decreasing property} holds for complex $V$-spaces, and Theorem \ref{main:isotrivial} follows from Theorem \ref{main} once we can prove that the moduli map $Y\to \mathfrak{M}_H$ associated to any smooth family $X\to Y\in\sM_H(Y)$ is a $V$-morphism.

For coarse moduli space of polarized Calabi-Yau manifolds, we can  prove a stronger result than Theorem \ref{main}.
\begin{main}\label{main:CY moduli}
Let $\mathfrak{M}_H$ be the coarse moduli space of polarized Calabi-Yau manifolds. Then there is  a K\"ahler $V$-metric on $\mathfrak{M}_H$ with non-positive holomorphic bisectional curvature and negative holomorphic sectional curvature. In particular, $\mathfrak{M}_H$ is Kobayashi $V$-hyperbolic.
\end{main}
For the definition of K\"ahler $V$-metric, we refer the reader to Definition \ref{def:V-metric}. By the Alhfors-Schwarz lemma, one can easily show that an orbifold equipped with a K\"ahler $V$-metric with negative holomorphic sectional curvature is Kobayashi $V$-hyperbolic.

Lastly, we prove a theorem which incorporates both the Shafarevich hyperbolicity conjecture and   Campana's isotriviality conjecture for smooth projective families of polarized Calabi-Yau manifolds, as suggested  by Campana  to the author.\noindent
\begin{main}\label{main:CY}
Let $f:X\to Y$ be a smooth projective family of   Calabi-Yau manifolds over a quasi-projective manifold $Y$.  
  Then the  \emph{essential dimension} ${\rm ess}(Y)$ is no less than   the \emph{variation of the family} ${\rm Var}(f)$, denoted to be the generic rank of the image of the Kodaira-Spencer map $T_Y\to R^1f_*T_{X/Y}$. In particular,
	\begin{enumerate}
		\item  \label{cyiso}if $Y$ is special, then $f$ is isotrivial. 
		\item \label{KK} Either $\kappa(Y)=-\infty$, or $\kappa(Y)\geq {\rm Var}(f)$. 
	\end{enumerate}  
\end{main}
We quickly show how to deduce the last two claims in the above theorem  from the main   result. By Campana \cite[Remarque 10.3.(1)]{Cam11}, the quasi-projective manifold is special if and only its essential dimension is $0$, which implies the vanishing of the Kodaira-Spencer map of $f$, and thus  \eqref{cyiso} follows.  By another theorem of Campana \cite[Example 10.4]{Cam11}, for any quasi-projective manifold $Y$, when $\kappa(Y)\geq 0$, one always has $\kappa(Y)\geq {\rm ess}(Y)$, which shows Theorem \ref{main:CY}.\eqref{KK}. Note that Theorem \ref{main:CY}.\eqref{KK} is a variant of a conjecture by Kebekus-Kov\'acs \cite[Conjecture 1.6]{KK08}. The original Kebekus-Kov\'acs  conjecture concerns  smooth families of canonically polarized manifolds, which was proved by Taji \cite{Taj16}.

Let us mention that, differently from the proof of Theorem \ref{main:CY moduli}, the proof of Theorem \ref{main:CY} does not use the work on coarse moduli spaces by Fujiki-Schumacher \cite{FS90}, and replies on the differential geometric properties of period mapping and period domains. In particular, we   give another proof of the isotriviality for smooth projective families of Calabi-Yau manifolds over {\it H}-special bases.
\subsection{Related results} 
In his remarkable paper on the birational classification of algebraic varieties, Campana introduced the \emph{special manifolds},  which are in a precise sense \emph{opposite} to manifolds of general type. 
\begin{definition}  A log pair $(X,D)$ is \emph{special} if  the strict inequality $\kappa(L)<p$ holds for all $p$ and all invertible sheaves  $L\subset \Omega^p_{X}(\log D)$. A quasi-projective manifold $Y$ is \emph{special} if there is a log pair $(X,D)$  such that $X$ is a smooth projective compactification of $Y$ with $D:=X-Y$ a  simple normal crossing divisor, and such that $(X,D)$ is special. 
\end{definition} 
	Important examples for special manifolds are rationally connected manifolds and projective manifolds $X$ with   Kodaira dimension $\kappa(X)=0$. Kobayashi has conjectured that   Calabi-Yau manifolds are all {\it H}-special, and this conjecture was generalized by Campana to all special manifolds \cite[Conjecture 9.2]{Cam04} in the same vein as the Lang conjecture.  
\begin{conjecture}[Campana]\label{conj:campana}
	A quasi-projective manifold $V$ is {special} if and only if it is {\it H}-special. 
\end{conjecture}
Conjecture \ref{conj:campana} is quite open at the present time, even for Calabi-Yau manifolds. It is only known for all $K3$ surfaces and   their Hilbert schemes by Kamenova-Lu-Verbitsky \cite{KLV14}.  

On the other hand, inspired by the Shafarevich-Viehweg-Zuo hyperbolicity conjecture, Campana made the following \emph{isotriviality conjecture}: any smooth projective family of canonically polarized manifolds over a \emph{special quasi-projective base manifold} is necessarily \emph{isotrivial}.   Jabbusch-Kebekus first proved this conjecture when the bases are surfaces \cite{JK11b} or threefolds \cite{JK11} by refining the Viehweg-Zuo sheaves in \cite{VZ02}. Combining this refined Viehweg-Zuo sheaves and Campana-P\u{a}un's work \cite{CP15,CP16,CP19},  Taji \cite{Taj16} completely  proved this conjecture. Therefore, Theorem \ref{main:isotrivial} provides new evidences for  Conjecture \ref{conj:campana}.

\medskip
Theorem \ref{main}  is inspired by the celebrated work of To-Yeung \cite{TY15}  on the Kobayashi hyperbolicity for  bases of effectively parametrized families of canonically polarized manifolds (see also \cite{BPW17,Sch18}). Indeed, as we will see in the proof of Theorem \ref{main}, we apply their methods of   \emph{augmented Weil-Petersson metric}  on these bases to prove a    \emph{Pick-Schwarz theorem} for   coarse moduli spaces of canonically polarized manifolds.

\medskip 

In \cite[Theorems A and B]{Den19}, we proved the following results which can be seen as an analytic version of the higher dimensional Shafarevich hyperbolicity conjecture.
\begin{theorem}\label{thm:Deng}
		Let $f:X\to Y$ be a smooth projective family  between   complex quasi-projective manifolds.
		\begin{enumerate}
			\item Assume that the general fiber  of $f$  has semi-ample canonical bundle, and $f$ is of maximal variation, namely the general fiber of $f$  can only be birational to  countably  other 
			fibers.  Then $Y$ is pseudo-Kobayashi hyperbolic.
			\item Assume that each fiber  of $f$  has big and nef  canonical bundle, and $f$ is effectively parametrized, namely the Kodaira-Spencer map is injective everywhere. 	Then  \(Y\) is Kobayashi hyperbolic.
		\end{enumerate}
	\end{theorem}
Therefore, inspired by Theorem \ref{main:isotrivial}, we propose the following conjecture.
\begin{conjecture}\label{conj}
	Let $f:X\to Y$ be a smooth projective family between   complex quasi-projective manifolds. Assume that each fiber  of $f$ has semi-ample canonical bundle, and the Kobayashi pseudo distance $Y$ vanishes identically. Then $f$ is birationally isotrivial, or even isotrivial. 
\end{conjecture}
Note that Conjecture \ref{conj} holds when $Y$ is a curve by Theorem \ref{thm:Deng}.
 \medskip
 
 Let us lastly mention an important work by  Borghesi-Tomassini \cite{BT17} which establishes the foundation of stack hyperbolicity. We think that it is worthwhile comparing  the $V$-hyperbolicity introduced in \S \ref{sec:V-hyper} with their work.
\subsection*{Conventions}
Throughout this paper we will work over complex number field $\bC$. A Calabi-Yau manifold is a simply connected complex projective
manifold with trivial canonical bundle. A fibration $f:X\to Y$ between compact K\"ahler manifolds is a surjective proper morphism with connected fibers. \emph{VPHS} is the abbreviation for \emph{variation of polarized Hodge structures}. A \emph{smooth projective family} of    quasi-projective manifolds is a smooth projective morphism with connected fibers.

\section{Preliminary}
 \subsection{Complex $V$-spaces}
By the work of Fujiki and Schumacher \cite{FS90}, \emph{complex $V$-space} arises naturally as the coarse moduli spaces of polarized manifolds. Let us  recall its precise definition.
\begin{definition}[\protect{\cite[Definition 1.4]{FS90}}]\label{def:V-space}
Let $X$ be a reduced complex space. Then
	\begin{enumerate}[leftmargin=0cm,itemindent=0.7cm,labelwidth=\itemindent,labelsep=0cm, align=left,label= {\rm (\roman*)},itemsep=0.07cm] 
	\item A \emph{local (analytic) Galois cover} of $X$ is a pair $\sF=(\pi:\tilde{U}\to U, G)$ consisting of a holomorphic map $\pi: \tilde{U}\to U$ of a connected complex space $\tilde{U}$ onto an open subset $U$ of
$X$ and of a finite group $G$ acting biholomorphically on $\tilde{U}$ over $U$ such that $\pi$ induces an
isomorphism $\tilde{U}/G\xrightarrow{\simeq}U$ of complex spaces.  
\item Let $\sF'=(\pi':\tilde{U}'\to U',G')$ be another local Galois cover of $X$ with $U'\subset U$, then
a morphism of $\sF'$  into $\sF$ is a commutative diagram of complex spaces
$$
\begin{tikzcd} 
\tilde{U}' \arrow[r, "\tilde{j}"] \arrow[d, "\pi'"]
& \tilde{U} \arrow[d, "\pi"] \\
U' \arrow[r, "j"]
& U 
\end{tikzcd}
$$
where $j$ is the inclusion and $\tilde{j}$ is an open embedding, together with an injective
homomorphism $\mu:G'\to G$ such that $\tilde{j}$ is $\mu$-equivariant. Any such $\tilde{j}$ is said to be associated to $j$. If such a morphism exists, $\sF$ and $\sF'$ are said to be compatible.
\item A \emph{complex $V$-structure} on $X$  is a collection
$$
\sV=\{\sF_\alpha=(\pi_\alpha:\tilde{U}_\alpha\to U_\alpha,G_\alpha)\}_{\alpha\in I}
$$
of local Galois covers of $X$ such that 
  $\{U_{\alpha}\}_{\alpha\in I}$  forms a \emph{basis} for open sets of  $X$,  and  
$\sF_\alpha$ and $\sF_\beta$  are compatible whenever $U_\alpha\subset U_\beta$.
\end{enumerate}
A \emph{complex $V$-space} is a complex space equipped with a $V$-structure. When each $\tilde{U}_\alpha$ is smooth, the definition of complex $V$-space reduces to the usual definition of  \emph{complex $V$-manifold} in the sense of Satake.
	\end{definition} 
Throughout the paper,  notations and terminologies in Definition \ref{def:V-space} will be used without mention once we deal with complex $V$-spaces.

 \begin{definition}[$V$-morphism]\label{def:V-morphism}
	Let $Y$ be an (ordinary) complex space, and let $X$ be a complex $V$-space. Then a morphism $f:Y\to X$ of complex spaces is a \emph{$V$-morphism} if it is \emph{locally liftable}. That is, for any point $y\in Y$, there exist a
neighborhood $\Omega$ of $y$ in $Y$, a local Galois cover $\pi:\tilde{U}\to U$ in $\sV$ with $f(z)\in U$, and a
morphism $\tau:\Omega\to \tilde{U}$ such that $f|_{\Omega}=\pi\circ \tau$.
\end{definition}
\subsection{Special manifold  and essential dimension}
Special varieties   and essential dimensions are originally introduced by Campana in  \cite{Cam04,Cam11}  studying the bimeromorphic structure of quasi-Kähler manifolds $X$. For any (quasi-)K\"ahler manifold $X$, Campana constructed a unique fibration (up to bimeromorphic equivalence), so-called  \emph{core map}, characterized
by having special orbifold fibres and general type orbifold base. The \emph{essential dimension} of $X$ is the dimension of the orbifold base of the core map. Instead of restating these   definitions from \cite{Cam11}, we recall a crucial result concerning the essential dimension which is sufficient for our purpose.  
\begin{lemma}\label{lem:essential dimension}
	Let $f:X\to Y$ be a fibration between  compact K\"ahler manifolds.  Assume that $D$ and  $\Delta$ are simple normal crossing divisors in $X$ and $Y$ so that $(f^*\Delta)_{\rm red}=D$, and $K_Y+\Delta$ is big. Then the essential dimension ${\rm ess}(X\setminus D)={\rm ess}(X,D)\geq {\rm dim}(Y)$.
\end{lemma}
 While Lemma \ref{lem:essential dimension}  is not explicitly stated in \cite{Cam11}, it can be easily obtained by combining    \cite[Theorems 3.5 and 5.3]{Cam11}, using \cite[Remark 3.8 and example 3.7]{Cam11}.
 
\section{$V$-hyperbolicity for  complex $V$-spaces}
\subsection{$V$-hyperbolicity}\label{sec:V-hyper}
In this section we introduce   \emph{$V$-hyperbolicity} for complex $V$-space, which follows the definition of \emph{hyperbolicity for $V$-manifolds} by Kobayashi \cite{Kob05} verbatim.
\begin{definition}
A complex $V$-space $X$ is \emph{Brody $V$-hyperbolic} if there exists no non-constant $V$-morphism $f:\bC\to X$. 
\end{definition}
To define the Kobayashi $V$-hyperbolicity for complex $V$-space, one has to introduce the Kobayashi pseudo distance for complex $V$-space. Recall that in the definition of the Kobayashi pseudo distance $d_X:X\times X\to \bR_{\geq 0}$ of a complex space $X$, one uses the chain of holomorphic mapping from the unit disk  $\bD$ to $X$. When $X$ is equipped with a $V$-structure, we only use the   $V$-morphism from the unit disk to $X$ to define the Kobayashi pseudo $V$-distance, denoted by $d^V_X:X\times X\to \bR_{\geq 0}$. Precisely speaking, for any $p,q\in X$,   a \emph{chain of $V$-disks} $\alpha$ connecting $p$ to $q$ is a set of $V$-morphisms  $\{f_i:\mathbb{D}\to X\}_{i=1,\ldots,m}$  so that $f_i(a_i)=p_{i-1}$, $f_i(b_i)=p_{i}$ for $\{a_i,b_i\}_{i=1,\ldots,m}$   pair of points in $\bD$ and  $p=p_0,q=p_m$. The length of $\alpha$ is defined to be  $\ell(\alpha)=\sum_{i=1}^{m}\rho(a_i,b_i)$, where $\rho:\bD\times\bD\to \mathbb{R}_{\geq 0}$ is the Poincar\'e metric of $\bD$. The Kobayashi pseudo $V$-distance from $p$ to $q$ is defined to be $d_X^V(p,q):=\inf_{\alpha}\ell(\alpha)$ where the infimum is taken over all chains of $V$-disks $\alpha$ connecting $p$ to $q$. It is easy to see that $d^V_X(p,q)\geq d_X(p,q)$ for any pair of points $p,q\in X$.
\begin{definition}
	A complex $V$-space $X$ is \emph{Kobayashi $V$-hyperbolic} if the Kobayashi pseudo $V$-distance  $d^V_X(p,q)>0$ for any pair of distinct points $p,q\in X$.
\end{definition}
It is easy to verify the following.
\begin{lemma}
	A complex $V$-space is Brody $V$-hyperbolic if it is Kobayashi $V$-hyperbolic.
\end{lemma}
The distance decreasing property holds for $V$-morphisms.
\begin{lemma}[Distance decreasing]\label{lem:decreasing}
	Let $X$ be a complex $V$-space, and let $f:Y\to X$ be a $V$-morphism from a reduced complex space $Y$.  Then $d_X^V(f(p),f(q))\leq d_Y(p,q)$  for any pair of distinct points $p,q\in Y$.
\end{lemma}
\begin{proof}
For	any chain of disks $\alpha$ connecting $p$ to $q$, its image under $f$, is a chain of $V$-disks from $f(p)$ to $f(q)$ by Definition \ref{def:V-morphism}. The lemma then follows from the definition of Kobayashi pseudo $V$-distance.
\end{proof}
\begin{lemma}\label{lem:largest}
For any complex $V$-space	$X$, the Kobayashi pseudo $V$-distance is the largest pseudo distance on $X$ satisfying the \emph{Schwarz-Pick property}. That is, for any pseudo distance $\delta:X\times X\to \bR_{\geq 0}$ such that
	\begin{eqnarray}\label{eq:decreasing}
	\delta_{X}(f(x),f(y))\leq \rho (x,y) 
	\end{eqnarray}
	for    $V$-disks $f:\bD\to X$, where $\rho:\bD\times \bD\to \bR_{\geq 0}$ is the Poincar\'e metric of the unit disk. Then $\delta_X\leq d_X^V$.
\end{lemma}
\begin{proof}
	Let $\alpha=\{f_i:\bD\to X\}_{i=1,\ldots,m}$  be a chain of $V$-discs from
	$p$ to $q$ so that $f_i(a_i)=p_{i-1}$, $f_i(b_i)=p_{i}$ for $\{a_i,b_i\}_{i=1,\ldots,m}$   pair of points in $\bD$ and  $p=p_0,q=p_m$. Then
	$$
	\delta_X(p,q)\leq \sum_{i=1}^{m}\delta_X(p_{i-1},p_{i})=\sum_{i=1}^{m}\delta_X(f_i(a_i),f_i(b_i)) \leq \sum_{i=1}^{m}\rho(a_i,b_i)=\ell(\alpha)
	$$
	where the first inequality is due to the triangle property of pseudo distance, and the third one follows from \eqref{eq:decreasing}. Hence
	$$
		\delta_X(p,q)\leq \inf_{\alpha} \ell(\alpha)=d_X^V(p,q).
	$$ 
\end{proof}
\subsection{Differential geometric criteria}
For a complex space $X$, its Zariski tangent space $T_X$ can be seen as a reduced complex space so that the natural projection map $\pi:T_X\to X$ is holomorphic.   A \emph{Finsler metric} on $X$ is a upper semi-continuous function $F:T_X\to \mathbb{R}_{\geq 0}$. Note that for
any point $x\in X$ there exist a closed embedding of a neighborhood $U$ of $x$ into a domain
$D\subset\bC^m$. Then $T_U\subset T_D$ is can be seen as a closed subvariety of the complex manifold $T_D$. In this setting,  
	A  Finsler (resp. hermitian) metric $h$ on a reduced complex space is  a \emph{$\sC^\infty$-metric} if there is a  $\sC^\infty$ Finsler (resp. hermitian) metric $h'$ of $T_D$, so that $h'|_{T_U}=h$.   
\begin{definition}\label{def:V-metric} 
	\begin{enumerate}[leftmargin=0cm,itemindent=0.7cm,labelwidth=\itemindent,labelsep=0cm, align=left,label= {\rm (\roman*)},itemsep=0.07cm] 
		\item 
 A Finsler (resp. hermitian, K\"ahler) $V$-metric on a  complex $V$-space is a collection $h=\{h_\alpha\}_{\alpha\in I}$ of $G_\alpha$-invariant  Finsler (resp. hermitian, K\"ahler)  metrics $h_{\alpha}$ on $\tilde{U}_{\alpha}$ such that if $U_\alpha\subset U_\beta$ then $\tilde{j}_{\alpha\beta}^*h_\beta=h_\alpha$  where $\tilde{j}_{\alpha\beta}:\tilde{U}_{\alpha}\to \tilde{U}_\beta$ is  
associated to the inclusion $U_\alpha\subset U_\beta$
$$
\begin{tikzcd} 
\tilde{U}_\alpha\arrow[r, "\tilde{j}_{\alpha\beta}"] \arrow[d, "\pi_\alpha"]
& \tilde{U}_\beta \arrow[d, "\pi_\beta"] \\
U_\alpha \arrow[r, "j_{\alpha\beta}"]
& U_\beta
\end{tikzcd}
$$
Such a  $V$-metric $h$ is \emph{smooth} if each $h_\alpha$ is a $\sC^\infty$  metric on $U_\alpha$.  
\item A K\"ahler $V$-metric $h$ has \emph{non-negative holomorphic bisectional curvature} if the holomorphic bisectional curvature of $h_\alpha$  is non-negative for each $\alpha\in I$. This K\"ahler $V$-metric $h$ has \emph{negative holomorphic  sectional curvature} if there is a constant $c>0$ so that the holomorphic sectional curvature of $h_\alpha$ is bounded from above by $-c$ for each $\alpha\in I$.
\end{enumerate}
\end{definition}
\begin{remark}
We stress that in Definition \ref{def:V-metric} hermitian and K\"ahler $V$-metrics are always required to be positively definite over each $\tilde{U}_\alpha$, while Finsler $V$-metric might be degenerate.
\end{remark}

Let $X$ be a complex $V$-space  equipped with a $\sC^\infty$ Finsler $V$-metric $h$. One can define a natural pseudo distance $d_{X,h}:X\times X\to \bR_{\geq 0}$ associated to $h$ as follows. For any two point $p,q\in X$, we set
$$
d_{X,h}(p,q):=\inf_{\gamma}\int_{0}^{1}|\gamma'(t)|_hdt
$$
where the infimum is taken over all piecewise smooth \emph{locally liftable} real curves $\gamma:[0,1]\to X$ joining $p$ and $q$. When $h$ is positively definite everywhere, $d_{X,h}$ is a distance.
\begin{lemma}\label{lem:basis}
	Let $X$ be a complex $V$-space and let $f:Y\to X$ be a $V$-morphism from a complex manifold $Y$ to $X$.	Then there exists  a \emph{basis} $\{Y_\beta\}_{\beta\in J\subset I}$ for open sets of $Y$ so that one has
	$$
	\begin{tikzcd}
	& \tilde{U}_\beta \arrow[dd,"\pi_\beta"]\\
Y_\beta \arrow[ur,"f_\beta"]\arrow[dr,"f|_{Y_\beta}"']&\\
	&  {U}_\beta
	\end{tikzcd}
	$$
	and if $Y_\gamma\subset Y_\beta$, there exists an element $g_\beta\in G_\beta$ with $\tilde{j}_{\gamma\beta}\circ f_\gamma=g_\beta\circ f_\beta|_{Y_\gamma}$.
\end{lemma}
 \begin{proof}
 By Definition \ref{def:V-space}, there is a collection
 $$
 \sV=\{\sF_\alpha=(\pi_\alpha:\tilde{U}_\alpha\to U_\alpha,G_\alpha)\}_{\alpha\in I}
 $$
 of local Galois covers of $X$ such that 
 $\{U_{\alpha}\}_{\alpha\in I}$  forms a \emph{basis} for open sets of  $X$,  and  
 $\sF_\alpha$ and $\sF_\beta$  are compatible whenever $U_\alpha\subset U_\beta$.	By Definition \ref{def:V-morphism}, there exists a collection of open sets $\{Y_\beta\}_{\beta\in I'\subset I}$ of $X$ so that $f|_{Y_\beta}:Y_\beta\to U_\beta$ can be lifted to a holomorphic map $ {f}_\beta:Y_\beta\to \tilde{U}_\beta$  so that $\pi_\beta\circ f_\beta=f|_{Y_\beta}$. 
 	Now we show that this collection can be refined further to be a basis. For any $\alpha,\beta\in I'$ with $Y_\alpha\cap Y_\beta\neq\emptyset$, one has $U_\alpha\cap U_\beta\neq\emptyset$. Pick any point $y\in Y_\alpha\cap Y_\beta$. Since $\{U_{\alpha}\}_{\alpha\in I}$  forms a \emph{basis} for open sets of  $X$, there exists $U_\gamma\subset U_\alpha\cap U_\beta$ with $\gamma\in I$ containing $f(y)$. 
 Pick a  connected open set $Y_\gamma\subset Y_\alpha\cap Y_\beta$ containing $y$	so that $f(Y_\gamma)\subset U_\gamma$. For $f_\beta(y)\in \tilde{U}_\beta$, it might happen that $f_\beta(y)\notin \tilde{j}_{\gamma\beta}(\tilde{U}_\gamma)$. However, one can find  an element $g_\beta\in G_\beta$ so that $g_\beta\cdot    f_\beta(y)\in \tilde{j}_{\gamma\beta}(\tilde{U}_\gamma)$. Since $\tilde{j}_{\gamma\beta}:\tilde{U}_\gamma\to \tilde{U}_\beta$ is an open embedding, we can shrink $Y_\gamma$ so that $g_\beta\circ  f_\beta(Y_\gamma)\subset \tilde{U}_\gamma$, where we think of $g_\beta:\tilde{U}_\beta\to \tilde{U}_\beta$ as a biholomorphism fo $\tilde{U}_\beta$. Note that $g_\beta\circ f_\beta:Y_\beta\to \tilde{U}_\beta$ is also a lift of $f|_{Y_\beta}$, that is, $\pi_\beta\circ g_\beta\circ f_\beta=f|_{Y_\beta}$. Write $f_\gamma=\tilde{j}_{\gamma\beta}^{-1}\circ(g_\beta\circ f_\beta)|_{Y_\gamma}:Y_\gamma\to \tilde{U}_\gamma$.  Hence $\pi_\gamma\circ f_\gamma=f|_{Y_\gamma}$, and one has $\tilde{j}_{\gamma\beta}\circ f_\gamma=g_\beta\circ f_\beta|_{Y_\gamma}$. 
 
 We have to show that there exists $g_{\alpha}\in G_\alpha$ so that $\tilde{j}_{\gamma\alpha}\circ f_\gamma=g_\alpha\circ f_\alpha|_{Y_\gamma}$. For any $g\in G_\alpha$, set $V_g:=\{y\in Y_\gamma\mid \tilde{j}_{\gamma\alpha}\circ f_\gamma(y)=g \circ f_\alpha (y) \}$, which is a   closed subvariety of $Y_\gamma$. Note that $$\pi_\alpha\circ\tilde{j}_{\gamma\alpha}\circ f_\gamma=\pi_\alpha \circ g \circ f_\alpha|_{Y_\gamma}=f|_{Y_\gamma}.$$
Hence for any $y\in Y_\gamma$, there exists   $g_y\in G_\alpha$ so that $y\in V_{g_y}$. In other words, $Y_\gamma=\cup_{g\in G_\alpha}V_g$. Since $Y_\gamma$ is a connected complex manifold, there exists $g_\alpha\in G_\alpha$ so that $Y_{\gamma}=V_{g_\alpha}$. One thus has $\tilde{j}_{\gamma\alpha}\circ f_\gamma=g_\alpha\circ f_\alpha|_{Y_\gamma}$. This accomplishes the proof.
 \end{proof}
The above lemma enables us to pull-back the $V$-metric via the $V$-morphism.
\begin{lemma}\label{lem:pull-back}
	Let $X$ be a complex $V$-space equipped with a $\sC^\infty$ Finsler $V$-metric $h$.	For any $V$-morphism $f:Y\to X$ from a complex manifold $Y$, one can pull-back the metric $h$ to $Y$, denoted by $f^*h$. 
\end{lemma}
\begin{proof}
Let  $\{Y_\beta\}_{\beta\in J\subset I}$  be a \emph{basis} for open sets of $Y$ as in  Lemma \ref{lem:basis}. Hence there exists a lift $f_\beta:Y_\beta\to \tilde{U}_\beta$ so that $\pi_\beta\circ f_\beta=f|_{U_\beta}$. By Definition \ref{def:V-metric}, over each $\tilde{U}_\beta$ there is a $G_\beta$-invariant smooth Finsler metric $h_\beta$ such that if $U_\alpha\subset U_\beta$ then $\tilde{j}_{\alpha\beta}^*h_\beta=h_\alpha$  where $\tilde{j}_{\alpha\beta}:\tilde{U}_{\alpha}\to \tilde{U}_\beta$ is  
associated to the inclusion $U_\alpha\subset U_\beta$. The pull-back $f_\beta^*h_\beta$ is a smooth Finsler metric on $Y_\beta$. We will show that $\{f_\beta^*h_\beta\}_{\beta\in J}$  can be glued together to form a smooth Finsler metric on $Y$. 

For any point $y\in Y_\alpha\cap Y_\beta$, there is  $\gamma\in J$ so that $Y_\gamma$ contains $y$ and $Y_\gamma\subset Y_\alpha\cap Y_\beta$. By Lemma \ref{lem:basis} one has  $\tilde{j}_{\gamma\beta}\circ f_\gamma=g_\beta\circ f_\beta|_{Y_\gamma}$. Hence $$f_\gamma^*h_\gamma=(\tilde{j}_{\gamma\beta}\circ  f_\gamma)^*h_\beta=(g_\beta\circ f_\beta)^*h_\beta|_{Y_\gamma}=f_\beta^*h_\beta|_{Y_\gamma}.$$
Here the first and the last inequality is due to Definition \ref{def:V-metric}. Therefore,  $\{f_\beta^*h_\beta\}_{\beta\in J}$ are compatible and they give rise to a metric on $Y$, denoted by $f^*h$.
\end{proof}

The following  criteria for Kobayashi $V$-hyperbolicity of complex $V$-spaces is important in the proof of Theorem \ref{main}.
\begin{proposition}\label{thm:infinitesimal}
Let $X$ be a complex $V$-space equipped with a $\sC^\infty$ Finsler $V$-metric $h$. Assume that for any $V$-morphism $f:\bD\to X$, one has $f^*h\leq h_{\bD}$, where $h_{\bD}$ is the infinitesimal Poincar\'e metric for the unit disk. Then $d_{X,h}\leq d_X^V$. If we further assume that $h$ is moreover positively definite, then $X$ is Kobayashi $V$-hyperbolic.
\end{proposition}
\begin{proof}
	By Lemma \ref{lem:largest}, it suffices to prove the Schwarz-Pick property for $d_{X,h}$. Let $f:\bD\to X$ be any $V$-morphism. For any $a,b\in \bD$, there exists a geodesic $\ell:[0,1]\in \bD$, so that  the Poinca\'e distance of $a$ and $b$ 
	$$
\rho(a,b)=\int_{0}^{1}|\ell'(t)|_{h_\bD}dt\geq  \int_{0}^{1}|\ell'(t)|_{f^*h}dt=\int_{0}^{1}|(f\circ\ell)'(t)|_{h}dt
	$$ 
	Since $f$ is a $V$-disk, $f\circ \ell:[0,1]\to X$ is a locally liftable real curve from $f(a)$ to $f(b)$, and thus $$\int_{0}^{1}|(f\circ\ell)'(t)|_{h}dt\geq d_{X,h}(f(a),f(b)).$$ Together with  Lemma \ref{lem:largest}, we prove the theorem.
\end{proof}

\section{Kobayashi $V$-hyperbolicity for coarse moduli space}
We first consider the (analytic) moduli functor $\sM_H$ of 
canonically polarized manifolds of dimension $n$ with the Hilbert polynomial $H(m)$, defined by
\begin{align*}
\sM_H(S):=\{&f: X\to S\mid f \mbox{ is proper smooth morphism with connected}\\ &\mbox{fibers};  K_F \mbox{ is ample  for each fiber }F, \chi(K_F^m)=H(m)\}/\simeq 
\end{align*}
where $``\simeq"$	is defined as follows: $(f_1:X_1\to S)\simeq (f_1:X_1\to S)$ if and only if there exists an $S$-isomorphism $X_1/S\stackrel{\simeq}{\to} X_2/S$. We also consider the   moduli functor $\sM_H$ of 
polarized Calabi-Yau manifolds of dimension $n$ with the Hilbert polynomial $H(m)$, defined by
\begin{align*}
\sM_H(S):=\{&(f: X\to S,\sL)\mid f \mbox{ is proper smooth morphism with} \\ & \mbox{simply connected fibers};  K_{X_s}=\sO_{X_s}; \sL \mbox{ is a line bundle on } X \\ & \mbox{ so that } L_s:=\sL|_{X_s}  \mbox{ is  ample},  \chi(L_s^m)=H(m)\}/\simeq 
\end{align*}
where $``\simeq"$	is defined as follows: $(f_1:X_1\to S,\sL_1)\simeq (f_1:X_1\to S,\sL_2)$ if and only if there exists an $S$-isomorphism $\phi:X_1/S\stackrel{\simeq}{\to} X_2/S$ so that $\sL_1\sim_{S}\phi^*\sL_2$.   By Fujiki-Schumacher \cite{FS90}, their   coarse moduli space $\mathfrak{M}_H$ exists and is a \emph{Hausdorff, reduced complex space}
with at most countably many connected components, carrying a natural $V$-structure.  Moreover,  coarse moduli space $\mathfrak{M}_H$ of polarized Calabi-Yau manifolds   
  is furthermore a complex $V$-manifold. Indeed,   any $[(X,L)]\in \mathfrak{M}_H$ has
a neighborhood in $\mathfrak{M}_H$ isomorphic to ${\rm Def}(X)/{\rm Aut}(X,L)$, where ${\rm Def}(X)$ is the Kuranishi
space of $X$, and ${\rm Aut}(X,L):=\{f:X\to X\mid f \mbox{ is biholomorphism}, f^*L=L \}$. By the work of Tian-Todorov, ${\rm Def}(X)$ is smooth. It is also well-known that ${\rm Aut}(X,L)$ is finite. Hence $\mathfrak{M}_H$ is a complex $V$-manifold. Let us mention that those coarse moduli spaces $\mathfrak{M}_H$ are moreover quasi-projective schemes by the deep work of Viehweg \cite{Vie95}, though we will not need this fact throughout this paper for we work in the category of analytic varieties and maps.
 
Theorem \ref{main} follows immediately if we can prove the \emph{infinitesimal Pick-Schwarz property}  for $\mathfrak{M}_H$ in Proposition \ref{thm:infinitesimal}. We will apply the deep work by   To-Yeung \cite{TY15} and Schumacher \cite{Sch12,Sch18}   to achieve this goal.  

For any  effectively parametrized smooth family  of canonically polarized manifolds,  To-Yeung \cite{TY15,TY18}   provides a \emph{canonical} way to define a Finsler metric (the so-called \emph{augmented Weil-Petersson metric} in \cite{TY18}) on the base with the holomorphic sectional curvature bounded away from zero, along the strategies initiated by Siu \cite{Siu86}. Let us first  review their construction. 

Consider a smooth proper morphism of canonically polarized manifolds $f:\sX\to S$ over the Riemann surface. 
In \cite{Siu86}  Siu defines a pseudo-hermitian metric  $h_1$ (so-called \emph{Weil-Petersson metric}) on $S$ using the $L^2$-inner product of the harmonic representative of the Kodaira-Spencer class with respect to the fiber-wise K\"ahler-Einstein metric. In general, the Gaussian curvature of $h_1$ is not negative. 
  The idea by  To-Yeung and Schumacher is to consider higher order Kodaira-Spencer map 
\begin{align}\label{eq:iterate}
\tau_p:T_S^{\otimes p}\to R^pf_*(\wedge^pT_{\sX/S})
\end{align}
to construct new metrics over $T_S$ by pulling back the  $L^2$-metric over $R^pf_*(\wedge^pT_{\sX/S})$.  Let $r>0$ be the largest integer so that $\tau_p$ is non-zero, which is called the \emph{length of the iterated Kodaira-spencer map}. Let $S_0\subset S$ be the complement of a discrete set in $S$ so that for each $p=1,\ldots, r$, $\dim H^p(\sX_s,\wedge^pT_{\sX_s})$ is constant  for any $s\in S_0$. Hence by the Grauert's direct image theorem, $R^pf_*(\wedge^pT_{\sX/S})|_{S_0}$ is locally free and carries a natural $L^2$ metric induced by the fiber-wise K\"ahler-Einstein metric. Pulling back this metric via $\tau_p$ and taking the $p$-th root, one obtains a $\sC^\infty$ pseudo hermitian metric over $S_0$, denoted by $h_p$ and called \emph{$p$-th order Weil-Petersson  metric}.   By the work of  Schumacher \cite{Sch12} and To-Yeung \cite{TY15}, these  metrics are \emph{canonical} in the following sense.
\begin{lemma}\label{lem:canonical}
For two smooth family of canonically polarized manifolds  $\{f_i:\sX_i\to S_i\}_{i=1,2}$ over the Riemann surface $S_i$ so that one has
$$
\begin{tikzcd}
\sX_1 \arrow[r,"F","\simeq"'] \arrow[d,"f_1"] & \sX_2  \arrow[d,"f_2"]\\
S_1 \arrow[r,"\mu","\simeq"'] & S_2 
\end{tikzcd}
$$
where both $F$ and $\mu$ are biholomorphisms. Then for the $p$-th order Weil-Petersson metrics $h_{i,p}$  over  $S_{i}'$   induced by $f_i$, one has $\mu^*h_{2,p}=h_{1,p}$. Here $S_{i}'\subset S_i$  is   the complement of a discrete set in $S_i$   so that for each $p=1,\ldots, r$, $\dim H^p(\sX_{i,s},\wedge^pT_{\sX_{i,s}})$ is constant  over $S'_i$. 
\end{lemma}  
Moreover, in \cite{TY15} To-Yeung  proved the following crucial  curvature formula  for $p$-th order Weil-Petersson metrics (see also Berntsson-P\u{a}un-Wang \cite{BPW17} and Schumacher \cite{Sch18}).
\begin{theorem}Let $h_p$ be the $p$-th order Weil-Petersson metric of $S_0\subset S$ associated to the smooth family $f:\sX\to S$ of canonically polarized manifolds of dimension $n$ over the Riemann surface $S$. Let $r>0$ be the largest integer so that $\tau_p$ in \eqref{eq:iterate} is non-zero. Set $K_p$ to be the Gaussian curvature of $h_p$. Then
\begin{align}
K_1&\leq -\frac{C_n}{c_1(K_{\sX_s})^n}+\frac{h_2}{h_1}\\
K_p&\leq \frac{1}{p}\big(-(\frac{h_p}{h_{p-1}})^{{p-1}}+(\frac{h_{p+1}}{h_{p}})^{{p+1}}\big) \mbox{ for } 1<p\leq r.
\end{align}
Here  $C_n$ is a positive constant depending only on $n$, and we make the convention that $h_p\equiv 0$ when $p>r$.
\end{theorem}
 Note that  $c_1(K_X)^n:X\in \mathfrak{M}_H\to \bZ_{>0}$ is a positive  constant, and  $c_1(K_{\sX_s})^n$ is thus independent of $s$.  Although none of the $p$-th order Weil-Petersson metric has negative Gaussian curvature,  To-Yeung \cite{TY15} established an algorithm to construct negative Gaussian curvature over $S_0$ by combining these Weil-Petersson metrics (see also \cite{BPW17,Sch12,Sch18}). That is, they prove that there exist constants $c_1,\ldots,c_r>0$ depending only on $n$, $c_1(K_{\sX_s})^n$ and $r$ so that the Gaussian curvature of the pseudo hermitian metric
\begin{eqnarray}\label{eq:combination}
h:=\sum_{i=1}^{r}c_ih_i
\end{eqnarray}
is bounded above by a negative constant $-\kappa$, where $\kappa$ also  only depends on $n$, $c_1(K_{\sX_s})^n$ and $r$. Let us now apply this deep result by To-Yeung   to prove Theorem \ref{main}.
\begin{proof}[Proof of Theorem \ref{main}]
 Let $\gamma:\bD\to \kM_H$ be any (non-constant) $V$-morphism. By Lemma \ref{lem:basis}, there exists a basis $\{C_{\beta}\}_{\beta\in J\subset I}$ for open sets   of $\bD$ so that  one has
$$
\begin{tikzcd}
& \tilde{U}_\beta \arrow[dd,"\pi_\beta"]\\
C_\beta \arrow[ur,"\gamma_\beta"]\arrow[dr,"\gamma|_{C_\beta}"']&\\
&  {U}_\beta
\end{tikzcd}
$$ 
For each $U_\beta$, there exists a  Kuranishi family $\tilde{f}_\beta:\tilde{\sX}_\beta\to \tilde{U}_\beta$ associated to certain canonically polarized manifold  $X_\beta\in \kM_H$. Any $g_\beta$ is induced by an automorphism of $\tilde{f}_\beta$. That is, there exists a biholomorphism $\tilde{g}_\beta:\tilde{\sX}_\beta\to \tilde{\sX}_\beta$ so that one has the following Cartesian square
\begin{equation}\label{eq:Cartesian}
\begin{tikzcd} 
\tilde{\sX}_\beta\arrow[r, "\tilde{g}_{\beta}", "\simeq"'] \arrow[d, "\tilde{f}_\beta"]
& \tilde{\sX}_\beta \arrow[d, "\tilde{f}_\beta"] \\
\tilde{U}_\beta \arrow[r, "g_{\beta}", "\simeq"']
& \tilde{U}_\beta
\end{tikzcd} 
\end{equation}
Moreover, by \cite[Lemma 5.5]{FS90}, if ${U}_\alpha\subset {U}_\beta$, then there is a morphism 
\begin{equation}\label{eq:non-Cartesian}
\begin{tikzcd} 
\tilde{\sX}_\alpha\arrow[r, "\varphi_{\alpha\beta}"] \arrow[d, "\tilde{f}_\alpha"]
& \tilde{\sX}_\beta \arrow[d, "\tilde{f}_\beta"] \\
\tilde{U}_\alpha \arrow[r, "\tilde{j}_{\alpha\beta}"]
& \tilde{U}_\beta
\end{tikzcd} 
\end{equation}
so that 
$\varphi_{\alpha\beta}:\tilde{\sX}_\alpha\xrightarrow{\sim} \tilde{\sX}_\beta|_{\tilde{f}_\beta^{-1}(\tilde{U}_\alpha)} $ is an isomorphism over $\tilde{U}_\alpha$.

Set $\sX_\beta:=\tilde{\sX}_\beta\times_{\tilde{U}_\beta}C_\beta$ and the induced smooth family $f_\beta:\sX_\beta\to C_\beta\in \sM_H(C_\beta)$. We will show that for the higher order Kodaira-Spencer map $$
\tau_{\beta,p}:T_{C_\beta}^{\otimes p}\to R^pf_*(\wedge^pT_{\sX_\beta/C_\beta})
$$
the largest integer of $p$ with $\tau_{\beta,p}$ non-zero, denoted by $r_{\beta}$, does not depend on $\beta\in J$. 

We first show that $r_\beta=r_\alpha$ if $C_\alpha\subset C_\beta$. By Lemma \ref{lem:basis}, there exists $g_\beta\in G_\beta$ so that $\tilde{j}_{\alpha\beta}\circ \gamma_\alpha=g_\beta\circ \gamma_\beta|_{C_\alpha}$. Note that the fiber product   of  $\tilde{f}_\beta:\tilde{\sX}_\beta\to \tilde{U}_\beta$ and $g_\beta\circ\gamma_\beta:C_\beta\to \tilde{U}_\beta$, is equal to $\sX_{\beta}\to C_\beta$. Indeed, this fact follows immediately from the diagram in \eqref{eq:Cartesian} is Cartesian.  Therefore, by \eqref{eq:non-Cartesian}, the family $\sX_\beta\times_{C_\beta}C_\alpha\to C_\alpha$ is isomorphic to $f_\alpha:\sX_\alpha\to C_\alpha$.  
 This shows the equality  $r_\beta=r_\alpha$. Since $\{C_\beta\}_{\beta\in J}$ is a basis for $\bD$, we conclude the equality $r_\beta=r_\gamma$ once $C_\beta\cap C_\gamma\neq\emptyset$. Set $\Upsilon_i:=\{C_\beta \mid r_\beta=i  \}$. Note that $1\leq r_\beta\leq n$ for any $\beta\in J$. For $i\neq j$, any $C_\beta\in \Upsilon_i$ and $C_\gamma\in \Upsilon_j$ satisfy $C_i\cap C_j=\emptyset$. Set $\bD_i:=\cup_{C_\beta\in \Upsilon_i}C_\beta$, which are open set of $\bD$, and one has $\bD=\coprod_{i=1}^{n}\bD_i$. Since $\bD$ is connected, we immediately conclude that $\bD_r=\bD$ for some $r\in \{1,\ldots,n\}$, and $\Upsilon_i=\emptyset$ for $i\neq r$. This shows that $r_\beta=r$ for any $\beta\in J$.

Write $h_{\beta,k}$ for the pseudo hermitian over $C_\beta^\circ$ induced by $\tau_{\beta,k}$, where $C_\beta^\circ\subset C_\beta$ the complement of a discrete set in $C_\beta$ so that  for each $p=1,\ldots, r$, $\dim H^p(\sX_{\beta,s},\wedge^pT_{\sX_{\beta,s}})$ is constant for $s\in C_\beta^\circ$. Recall that the family $\sX_\beta\times_{C_\beta}C_\alpha\to C_\alpha$ is isomorphic to $f_\alpha:\sX_\alpha\to C_\alpha$. Then $C_\alpha^\circ=C_\beta^\circ\cap C_\alpha$ if $C_\alpha\subset C_\beta$, and one has $h_{\beta,k}|_{C_\alpha^\circ}=h_{\alpha,k}$ by Lemma \ref{lem:canonical}. 

 Recall that  $c_1(K_X)^n:X\in \mathfrak{M}_H\to \bZ_{>0}$ is a constant. Hence the coefficients  $c_1,\ldots,c_r$ in the algorithm constructing the pseudo hermitian metric $h_\beta:=\sum_{k=1}^{r}c_kh_{\beta,k}$ over $C_\beta^\circ$ in \eqref{eq:combination} with negative Gaussian curvature do  not depend on $\beta\in J$. That is, writing $h_\beta:=\sqrt{-1}\varphi_\beta(z)dz\wedge d\bar{z}$, where $\varphi_\beta(z)$ is a smooth non-negative function on $C_\beta^\circ$, one has 
 \begin{align}\label{eq:psh}
 \frac{\partial^2\log \varphi_\beta(z)}{\partial z {\partial}\bar{z}}\geq  \kappa_r\varphi_\beta(z)
 \end{align} for some positive constant $\kappa_r$ not depending on $\beta$. By \cite[Proposition 12]{Sch12}, $\varphi_\beta(z)$ is bounded from above on any relatively compact subsets in $C_\beta$, and thus it extends to a psh function over $C_\beta$, still denoted by $\varphi_\beta$ abusively, so that \eqref{eq:psh} holds over $C_\beta$. Recall that $h_{\beta,k}|_{C_\alpha^\circ}=h_{\alpha,k}$ for any $C_\alpha\subset C_\beta$, and $\{C_\beta\}_{\beta\in J}$ is a basis for open sets in $\bD$.   Hence  the psh functions $\{\varphi_\beta\}_{\beta\in J}$ can be glued together to obtain a psh function $\varphi\in {\rm Psh}(\bD)$, so that the Gaussian curvature of the pseudo hermitian metric $h:=\sqrt{-1}\varphi(z)dz\wedge d\bar{z}$ is bounded above by $-\kappa_r$.  By Demailly's Ahlfors-Schwarz lemma \cite{Dem97}, $h\leq \frac{1}{\kappa_r}h_\bD$. 
 
 By \cite[Theorem 7.10]{FS90}, there is a $\sC^\infty$ K\"ahler $V$-metric $h_{WP}$ on $\mathfrak{M}_H$ induced by $L^2$-inner product of the harmonic representative of the Kodaira-Spencer class with respect to the fiber-wise K\"ahler-Einstein metric.  It follows from Lemma \ref{lem:pull-back} that one can pull back $h_{\rm WP}$ via $\gamma$ to define a pseudo hermitian metric $\gamma^*h_{\rm WP}$  on $\bD$. By the construction of $h_{WP}$ in \cite[\S 7]{FS90}, one has $h_{\beta,1}=\gamma^*h_{WP}|_{C_\beta^\circ}$ for any $\beta\in J$.  Note that $\varphi(z)$ is upper semi-continuous, and $h|_{C_\beta^\circ} \geq h_{\beta,1}=\gamma^*h_{WP}|_{C_\beta^\circ}$ for any $\beta\in J$. Hence one has $h\geq \gamma^*h_{WP}$ and  $h_{\bD}\geq \gamma^*(\kappa_r\cdot h_{WP})$. Let us stress that the positive constant $\kappa_r$ only depends on the length $r$ of the iterated Kodaira-spencer map of \eqref{eq:iterate}.
 
 For other non-constant $V$-disk, the length $r$ of the iterated Kodaira-spencer map   might be different. However, they can only have value in $\{1,\ldots,n\}$. Set $\kappa:={\rm min}\{\kappa_1,\ldots,\kappa_r \}>0$. Then for any $V$-disk $\tilde{\gamma}:\bD\to \mathfrak{M}_H$, one always has $h_{\bD}\geq \gamma^*(\kappa\cdot h_{WP})$. The theorem follows from Proposition \ref{thm:infinitesimal} immediately. 
\end{proof}
For the coarse moduli space $\mathfrak{M}_H$ of polarized Calabi-Yau manifolds, we will show that the Hodge metric will induce a natural K\"ahler $V$-metric  $\mathfrak{M}_H$ with (semi-)negative holomorphic (bi)sectional curvature.
\begin{proof}[Proof of Theorem \ref{main:CY moduli}]   As is well-known, the coarse moduli space $\mathfrak{M}_H$ of polarized Calabi-Yau manifolds are complex $V$-manifolds. Indeed, for any polarized Calabi-Yau manifold $ (X,L) $,  its automorphism group   is finite, and the polarized Kuranishi spaces of $(X,L)$ are smooth by the work of Tian-Todorov.  We will use the Hodge metric to construct a $\sC^\infty$ \emph{K\"ahler $V$-metric} over $\mathfrak{M}_H$.

By the construction of $\mathfrak{M}_H$ in \cite{FS90}, each $\tilde{U}_\alpha$ is the base manifold of the polarized Kuranishi family $(f_\alpha:{\sX}_\alpha\to\tilde{U}_\alpha,\sL_\alpha)\in \sM_H(\tilde{U}_\alpha)$  of some  polarized  Calabi-Yau manifold $(X_\alpha,L_\alpha)$ and $G_\alpha={\rm Aut}(X_\alpha,L_\alpha)$ with $[(X_\alpha,L_\alpha)]\in \mathfrak{M}_H$. We can assume  each $\tilde{U}_\alpha$ is a simply connected complex manifold.  Fixing a marking 
$$m_\alpha:H_{\rm prim}^n(X_\alpha,\bZ)/{\rm torsion}\stackrel{\simeq}{\to}\bZ^r,$$
we can define a period map $p_\alpha:\tilde{U}_\alpha\to \Omega_\alpha$ associated to the  variation of polarized Hodge structure of this polarized Kuranishi family. Since the Kodaira-Spencer map $T_{\tilde{U}_\alpha}\to R^1(f_\alpha)_*T_{\sX_\alpha/\tilde{U}_\alpha}$ is an isomorphism,   $p_\alpha$ is an immersion by the infinitesimal Torelli theorem of Calabi-Yau manifolds. Moreover, $p_\alpha$ is horizontal, and one can pull back the Hodge metric of $T_{\Omega_\alpha}^{-1,1}$ via $p_\alpha$ to $\tilde{U}_\alpha$, denoted by $h_\alpha$.  By \cite{Lu99}, $h_\alpha$ a K\"ahler metric over $\tilde{U}_\alpha$.   For any  $g_\alpha\in G_\alpha$ acting biholomorphically on $\tilde{U}_\alpha$,  one has
the following commutative diagram 
\begin{equation}\label{eq:non-Cartesian3}
\begin{tikzcd} 
{\sX}_\alpha \arrow[r] \arrow[d, "{f}_\alpha"]
&  {\sX}_\alpha  \arrow[d, "{f}_\alpha"] \\
\tilde{U}_\alpha \arrow[r, "g_\alpha"]
& \tilde{U}_\alpha
\end{tikzcd} 
\end{equation}
which is  Cartesian. Hence one has $p_\alpha\circ g_\alpha=p_\alpha$ which in particular shows that $h_\alpha$ is $g_\alpha$ invariant.  We will show that those metrics $\{h_\alpha\}_{\alpha\in J}$ can be glued together to obtain a $V$-metric over $\mathfrak{M}_H$.
 
By \cite[Lemma 5.5]{FS90} again, if ${U}_\alpha\subset {U}_\beta$, then there is a morphism 
\begin{equation}\label{eq:non-Cartesian2}
\begin{tikzcd} 
 {\sX}_\alpha \arrow[r, "\varphi_{\alpha\beta}"] \arrow[d, "{f}_\alpha"]
&  {\sX}_\beta  \arrow[d, "{f}_\beta"] \\
\tilde{U}_\alpha \arrow[r, "\tilde{j}_{\alpha\beta}"]
& \tilde{U}_\beta
\end{tikzcd} 
\end{equation}
so that 
$\varphi_{\alpha\beta}: {\sX}_\alpha\xrightarrow{\sim}  {\sX}_\beta|_{{f}_\beta^{-1}(\tilde{U}_\alpha)} $ is an isomorphism over $\tilde{U}_\alpha$, with $\varphi_{\alpha\beta}^*\sL_\beta=\sL_\alpha$.  This implies that   $\Omega_\alpha=\Omega_\beta$ for any $U_\alpha\subset U_\beta$. Note that $\{U_\alpha\}_{\alpha\in J}$ forms a basis of open sets of $\mathfrak{M}_H$. Fixing any   connected component $\Gamma$ of $\mathfrak{M}_H$, one has $\Omega_\alpha=\Omega_\beta$ if $U_\alpha$ and $U_\beta$ are open sets in $\Gamma$. We thus write $\Omega$ for $\Omega_\alpha$ when no confusion can arise.

To show that $\{h_\alpha\}_{\alpha\in J}$ is a K\"ahler $V$-metric over $\mathfrak{M}_H$, one has to prove that $\tilde{j}_{\alpha\beta}^*h_\beta=h_\alpha$. This is not obvious for the reference point and the marking defining the period map $p_\alpha$ might not be compatible with those of $p_\beta$ via the inclusion $\tilde{j}_{\alpha\beta}$. It indeed follows from the following result.
\begin{claim}
	Let $(f_1:\sX_1\to S_1, \sL_1)$ and $(f_2:\sX_2\to S_2,\sL_2)$ be two smooth family of polarized manifolds of relative dimension $n$. Assume that there are biholomorphisms  $\phi:\sX_1\to\sX_2$ and $\varphi:S_1\to S_2$ so that one has the following commutative diagram
	\begin{equation}\label{eq:non-Cartesian4} 
	\begin{tikzcd} 
	{\sX_1}  \arrow[r, "\phi","\simeq"'] \arrow[d, "f_1"]
	&  {\sX_2}   \arrow[d, "{f_2}"] \\
S_1 \arrow[r, "\varphi","\simeq"']
	& S_2
	\end{tikzcd} 
	\end{equation}
and  $\phi^*\sL_2\equiv_{S_1} \sL_1$, where $\equiv_{S_1}$ stands for relatively numerically equivalent. Let $h_{S_i}$ be the Hodge metric  on $S_i$   induced by   period maps of the variation of polarized Hodge structures  $(\bV_i:=R^n_{\rm prim}(f_i)_*\bZ,F^\bullet\cV_i,Q_i)$. Here $\bV_i$ is the local system of primitive cohomology, $F^\bullet\cV_i$ is the Hodge filtration of $\cV_i:=\bV_i\otimes_{\bZ}\sO_{S_i}$, and  $Q_i:\bV_i\times \bV_i\to \bZ$ is the bilinear form induced by the polarization $\sL_i$.
Then   $(\varphi_*)^{-1}h_2=h_1$.  In particular, the Hodge metric does not depend on the reference point and the marking.
\end{claim}
\begin{proof}
We first observe that	for any isomorphism  between polarized manifolds $\psi:(X_1,A_1)\to (X_2,A_2)$ of dimension $n$ so that $\psi^*A_2\equiv A_1$,  the isomorphism $\psi^*:H_{\rm prim}^m(X_2,\bC)\stackrel{\simeq}{\to} H_{\rm prim}^m(X_1,\bC)$ is  an isometry  with respect to the Hodge metric. Indeed, for any $\alpha,\beta\in H^m(X_i,\mathbb{R})$, one defines 
\begin{align*} 
Q_i(\alpha,\beta) :=(-1)^{\frac{n(n-1)}{2}}\alpha\cdot\beta\cdot c_1(A_i)^{n-m}\in \mathbb{R} 
\end{align*}
and the Hodge metric $h_i$ over $H^m_{\rm prim}(X_i,\mathbb{R})$ is defined by $$h_i(\alpha,\beta):=\sum_{p=0}^{m}\sqrt{-1}^{2p-m}Q(\alpha_p,\overline{\beta_{p}}),$$ where $\alpha=\sum_{p=0}^{m}\alpha_p$ is the component of the Hodge composition of $\alpha$ with $\alpha_p\in H^{p,m-p}_{\rm prim}(X_i):=H^{p,m-p}(X_i)\cap H_{\rm prim}^m(X_i,\bC)$.
It is easy to verify that $Q_1(\psi^*\alpha,\psi^*\beta)=Q_2(\alpha,\beta)$.
Since $\psi^*:H^{m}(X_2,\bC)\stackrel{\simeq}{\to} H^{m}(X_1,\bC)$ keeps the Hodge decomposition and the primitive cohomology invariant, one thus has
$$
h_2(\alpha,\beta)=h_1(\psi^*\alpha,\psi^*\beta).
$$
Note that the biholomorphism $\phi:\sX_1\to \sX_2$ induces an isomorphism between $(\bV_i,F^\bullet\cV_i,Q_i)$. Hence for the corresponding Hodge bundles  $$\sH_i:=Gr^F_{\bullet} \cV_i=\oplus_{p=0}^{n}R_{\rm prim}^p(f_i)_*(\Omega^{n-p}_{\sX_i/S_i})$$ equipped with the Hodge metric $h_i$, one has $\phi^*:\sH_2\stackrel{\simeq}{\to} \sH_1$. By the above argument, one has $\phi^*h_2=h_1$.  For the natural Higgs structures $(\sH_i,\theta_i)$ with $\theta_i:\sH_i\to \sH_i\otimes \Omega_{S_i}$ and $\theta_i\wedge\theta_i=0$,  one has
	\begin{equation} \label{dia:induce}
\begin{tikzcd} 
\sH_2  \arrow[r, "\theta_2"] \arrow[d, "\phi^*"]
&  \sH_2\otimes \Omega_{S_2}  \arrow[d, "\phi^*\otimes \varphi^*"] \\
\sH_1 \arrow[r, "\theta_1"]
& \sH_1\otimes \Omega_{S_1}
\end{tikzcd} 
\end{equation} 
The Higgs bundle $(\sH_i,\theta_i)$ induces a morphism $\eta_i:T_{S_i} \to {\rm End}(\sH_i,\sH_i)$. Let $h_{{\rm End}(\sH_i)}$ be the smooth metric on ${\rm End}(\sH_i)$ induced by $h_i$. Then for the pseudo K\"ahler metric $h_{S_i}$ on $T_{S_i}$   induced by the pull-back of the metric $h_{{\rm End}(\sH_i)}$ via  $\eta_i$, one has $(\varphi_*)^{-1}h_{S_2}=h_{S_1}$ via   the isomorphism $\varphi_*:T_{S_1}\stackrel{\simeq}{\to} T_{S_2}$.
\end{proof}
The above claim shows that $h:=\{h_\alpha\}_{\alpha\in J}$ is a smooth K\"ahler $V$-metric on $\mathfrak{M}_H$.  By the work of Griffiths-Schmid \cite{GS69} and Peters \cite{Pet91} (see also \cite[Theorem 13.6.3]{CMP17}), the holomorphic sectional curvature of a period domain $\Omega$ is 
negative and uniformly bounded away from 0 in the horizontal direction $T_{\Omega}^{-1,1}$. By the curvature decreasing property, we conclude that the holomorphic bisectional curvature of $h_\alpha$ is non-positive for each $\alpha\in J$. Note that $\Omega_\alpha=\Omega_\beta$ for  any $U_\alpha$ and $U_\beta$ contained in  the same connected component $\Gamma$ of $\mathfrak{M}_H$. Hence there is a constant  $c>0$ so that the holomorphic sectional curvature of $h_\alpha$ is bounded from above by $-c$ for any $U_\alpha\subset \Gamma$. This proves the main result of the theorem. 

Let us prove the Kobayashi $V$-hyperbolicity of $\mathfrak{M}_H$. Let $\gamma:\bD\to \mathfrak{M}_H$ be  any $V$-disk. By Lemma \ref{lem:pull-back}, one can pull back the above K\"ahler $V$-metric on $\mathfrak{M}_H$, denoted by $\gamma^*h$, which is a smooth hermitian metric on $\bD$. Since the holomorphic sectional curvature of $h_\alpha$ is bounded above by $-c$ for each $\alpha\in J$, by the definition of $\gamma^*h$, the Gaussian curvature of $\gamma^*h$ is bounded above by $-c$. By the Ahlfors-Schwarz lemma, one has $\gamma^*(ch)\leq h_{\bD}$, where $h_{\bD}$ is the infinitesimal Poincar\'e metric on $\bD$. It follows from Proposition \ref{thm:infinitesimal} that $\mathfrak{M}_H$ is Kobayashi $V$-hyperbolic.
\end{proof}
Let us show how to conclude Theorem \ref{main:isotrivial} from Theorems \ref{main} and \ref{main:CY moduli}.
\begin{proof}[Proof of Theorem \ref{main:isotrivial}]
 Denote $n$ to be the relative dimension of the fibration $f$. Suppose by contradiction that $f$ is not isotrivial. We  prove the theorem for smooth families of canonically polarized manifolds, and the proof  for the case of smooth families of polarized Calabi-Yau manifolds is exactly the same. The Euler characteristic
 $$\chi(mK_F)=\sum_{i=0}^{n}(-1)^ih^i(X,mK_F)=h^0(X,mK_F)\mbox{ for } m\geq 2$$
 which is invariant under smooth deformation.  Hence there exists a polynomial $H(m)$ so that $\chi(mK_F)=H(m)$ for each fiber $F$. Consider the  moduli functor $\sM_H$ of  
 canonically polarized manifolds of dimension $n$ with the Hilbert polynomial $H(m)$, and by the work of Fujiki-Schumacher,  there exists a coarse moduli space $\mathfrak{M}_H$   of $\sM_H$, which is a Hausdorff complex $V$-space. By Theorem \ref{main}, $\mathfrak{M}_H$  is Kobayashi $V$-hyperbolic.  Then the moduli map $\varphi:Y\to \mathfrak{M}_H$ is not   constant. Since   $\mathfrak{M}_H$  is locally defined by quotient of Kuranishi spaces, by the universal property of Kuranishi families, $\varphi$ is locally liftable, and thus is a $V$-morphism. Pick a pair of points $p,q\in Y$ so that $f(p)\neq f(q)$. By Lemma \ref{lem:decreasing}, the Kobayashi pseudo distance $d_Y(p,q)\geq d_{\mathfrak{M}_H}^V(f(p),f(q))>0$, which contradicts to that $Y$ is {\it H}-special. 
\end{proof}

\section{Essential dimension for bases of  Calabi-Yau families}
The proof of Theorem \ref{main:CY} is inspired by the work \cite{Bru18,BC17,BBT18}.
\begin{proof}[Proof of Theorem \ref{main:CY}]
The theorem is trivial if ${\rm Var}(f)=0$. Assume now ${\rm Var}(f)>0$. 
 Let $n$ be the relative dimension of $f:X\to Y$. Since $f$ is a smooth projective family, one can thus   fix a polarization $(f:X\to Y,L)$.   Consider the $\mathbb{Z}$-variation  of polarized Hodge structure   $(\bV:=R^n_{\rm prim}(f)_*\bZ,F^\bullet\cV,Q)$,  where $\bV$ is the local system of primitive cohomology, $F^\bullet\cV$ is the Hodge filtration of $\cV :=\bV \otimes_{\bZ}\sO_{Y}$, and  $Q:\bV\times \bV\to \bZ$ is the bilinear form induced by the polarization $L$.  Since the infinitesimal Torelli theorem holds for Calabi-Yau manifolds, the period map $p:Y\to \Gamma\backslash \Omega$ of the variation of polarized Hodge structures $(\bV,F^\bullet\cV,Q)$ is not constant unless the Kodaira-Spencer map is identically zero. Here $\Omega$ denotes to be the polarized period domain, and $\Gamma$ is the monodromy group acting on $\Omega$.  As in \cite{Som78}, we can pass to a finite \'etale cover $e:{Y}' \to Y$, so that for the induced period map $p':Y'\to \Gamma'\backslash \Omega$ of $e^*(\bV,F^\bullet \cV,Q)$, the  monodromy $\Gamma'$ acting on $\Omega$  is torsion free. Since $\Gamma'$ acts on $\Omega$  discretely, it must be free. Take a smooth projective manifold   $\overline{Y'}$ which compactifies $Y'$ with $\overline{Y'}-Y'$   a simple normal crossing divisor. Then the local monodromy of $e^*\bV$ around the boundary is either trivial or infinite.  By a theorem of Griffiths \cite{Gri70},  the $\mathbb{Z}$-variation  of polarized Hodge structure  extends over the locus where the local monodromy is trivial, and we denote by $\tilde{Y}$ be the largest open set of $\overline{Y'}$ so that $e^*(\bV,F^\bullet \cV,Q)$ extends. Then  $\tilde{Y}$ is the Zariski open set of $\overline{Y}$, and the  period map of the extended VPHS is denoted to be $\tilde{p}:\tilde{Y}\to \Gamma'\backslash \Omega$.  By another theorem of Griffiths \cite{Gri70}, the period map  of VPHS    extends  as a proper map over the locus where the local monodromy is finite. Hence $\tilde{p}:\tilde{Y}\to \Gamma'\backslash \Omega$ is proper.  Let $\tilde{Y}\xrightarrow{\tilde{\pi}}\tilde{Z}\xrightarrow{\mu} \Gamma'\backslash \Omega$ be a Stein factorization of $\tilde{p}$. By \cite[Proposition \rom{3}, Remark \rom{3}.C]{Som78} (see also \cite[Theorem 7.7]{BBT18}), one has
	$$
	\begin{tikzcd}
	\tilde{Y} \arrow[d,"\tilde{\pi}"]&  Y^\circ \arrow[d,,"{\pi}^\circ"]\arrow[l,"\alpha"']\arrow[r,hook] & \overline{Y}\arrow[d,"\overline{\pi}"]\\
		\tilde{Z}  &  Z^\circ  \arrow[l,"\beta"']\arrow[r,hook] & \overline{Z}
	\end{tikzcd}
	$$
	where $\alpha$ and $\beta$ are proper modifications,   $Y^\circ$ and $Z^\circ$ are Zariski dense open sets of compact K\"ahler manifolds  $\overline{Y}$ and $\overline{Z}$ respectively, and $\pi^\circ:=\overline{\pi}|_{Y^\circ}$ is a surjective morphism with connected fibers. We can assume that $D=\overline{Y}-Y^\circ$ and $E=\overline{Z}-Z^\circ$ are both simple normal crossing divisors.  The image   of $j:=\mu\circ\beta $  is a closed analytic subvariety of $\Gamma'\backslash\Omega$. 
	Note that the composition $j:Z^\circ\to \Gamma'\backslash\Omega$ is always locally liftable to $\Omega$ for $\Gamma'$ acts freely over $\Omega$. Since $\tilde{p}$ is horizontal, the  map $j:Z^\circ\to \Gamma'\backslash \Omega$ is also horizontal. Indeed,  the embedding of the regular locus of $\tilde{p}(\tilde{Y})$ to $\Gamma'\backslash\Omega$ is horizontal for $\tilde{p}$ is horizontal, and by continuity $j:Z^\circ\to \Gamma'\backslash\Omega$ is horizontal as well. By \cite[Lemma-Definition 4.6.3]{CMP17}, $j:Z^\circ\to \Gamma'\backslash\Omega$ is the period map of another VPHS, which is generically immersive for  $\mu$ is finite. 
	
  By Campana \cite[Proposition 10.11]{Cam11}, the essential dimension of an orbifold is invariant under \'etale cover and bimeromorphic modification, and thus ${\rm ess}(Y')={\rm ess}(Y)$. It follows from $Y'\subset \tilde{Y}\subset \overline{Y'}$ that ${\rm ess}(Y^\circ)={\rm ess}(\tilde{Y})\leq {\rm ess}(Y')$.   
 Recall that $j:Z^\circ\to \Gamma'\backslash\Omega$ is the period map of a  VPHS, which is generically immersive.	By an elegant work of Cadorel-Brunebarbe \cite{BC17}, $Z^\circ$ is of log general type, or in other words, $K_{\overline{Z}}+E$ is big. Since $\overline{\pi}^{-1}(E)=D$,  
  it follows from Lemma \ref{lem:essential dimension} that  
  ${\rm ess}(Y^\circ)={\rm ess}(\overline{Y},D)\geq \dim \overline{Z}$. Since the infinitesimal Torelli theorem holds for Calabi-Yau manifolds, one thus has $\dim \overline{Z}=\dim \tilde{Z}={\rm Var}(f)$. This proves that ${\rm ess}(Y)\geq {\rm Var}(f)$.  

In the rest of the proof we study the pseudo Kobayashi distance of $Y$.   The $\Gamma'$-invariant natural metric $h$ on period domain $\Omega$ induces a pseudo K\"ahler metric (in particular, a $\sC^\infty$ Finsler metric) $h_{Z^\circ}$ over $Z^\circ$. By Griffiths-Schmid  (see  also \cite[\S 13.6]{CMP17}), 
	the holomorphic  sectional   curvature  of $h$ along the holomorphic horizontal tangent bundle $T_{\Omega}^{-1,1}$ is negative and uniformly bounded away from 0  by a negative constant. By curvature decreasing property, the holomorphic sectional curvature of $h_{Z^\circ}$ is also bounded from above by a negative constant. Let $\Delta\subset  Z^\circ$ be the locus where $j$ is not immersive, which is a closed subvariety of $Z^\circ$.  Then $h_{Z^\circ}$ is positively definite over the Zariski open set $U:=Z^\circ\setminus \Delta$.  By Ahlfors-Schwarz the Kobayashi pseudo distance $d_{Z^\circ}(p,q)>0$ for any pair of distinct points $p,q\in Z$ not both  contained in $\Delta$ (see \emph{e.g.} \cite[Theorem 3.7.4]{Kob98}). Since $\pi^\circ:Y^\circ\to Z^\circ$ is proper and surjective, one can find  two points $x,y\in Y^\circ$ such that $\pi^\circ(x)\neq \pi^\circ(y)$, and such that $\pi^\circ(x)$ and $\pi^\circ(y)$ are not both contained in $\Delta$. By distance decreasing property of Kobayashi pseudo distance, one has $d_{Y^\circ}(x,y)\geq d_{Z^\circ}\big(\pi^\circ(x),\pi^\circ(y)\big)>0$. 
	
Assume that  $Y$ is {\it H}-special. Since $Y'\to Y$ is a finite \'etale cover, $Y'$ is also    {\it H}-special. Recall that $Y'\subset \tilde{Y}\subset \overline{Y'}$. 
Then by the distance decreasing propery the Kobayashi pseudo distance $d_{\tilde{Y}}|_{Y'}\leq d_{Y'}\equiv 0$. By the continuity of Kobayashi pseudo distance \cite[Proposition 3.1.13]{Kob05} we conclude that $d_{\tilde{Y}}\equiv 0$. Since $\alpha:Y^\circ\to \tilde{Y}$ is proper modification, $d_{Y^\circ}\equiv 0$. We obtain the contradiction. 
\end{proof}
\medskip

\noindent \textbf{Acknowledgments.} I sincerely thank Professors Fr\'ed\'eric Campana, Jean-Pierre Demailly, Ken-Ichi Yoshikawa, and Yohan Brunebarbe, Junyan Cao, Philipp Naumann for answering my questions and very helpful discussions.   
	 

\end{document}